\documentclass[12pt]{article}
\usepackage{amssymb}

\usepackage{amsmath}
\usepackage{amsmath,amsxtra,amssymb,latexsym, amscd,amsthm}

\advance\voffset-1truecm\relax
\advance\hoffset-1truecm\relax
\def\C{{\Bbb C}}

\begin{document}

\title{UNIQUENESS PROBLEM\\
FOR MEROMORPHIC MAPPINGS\\
 WITH TRUNCATED MULTIPLICTIES\\
 AND FEW TARGETS}
\author{$\quad $Gerd Dethloff and Tran Van Tan}
\date{$\quad$}
\maketitle

\begin{abstract}
\noindent In this paper, using techniques of value distribution theory, we give a
uniqueness theorem for meromorphic
mappings of $\mathbb{C}^{m}$ into $\mathbb{C}P^{n}$ 
with truncated multiplicities and "few" targets. We also give 
a theorem of linear degeneration  for such maps with truncated
multiplicities and moving targets.\\
\begin{center}
{\bf R\'{e}sum\'{e}}
\end{center}
Dans cet article, on donne un th\'{e}or\`{e}me d'unicit\'{e} pour des applications
m\'{e}romorphes de $\mathbb{C}^{m}$ dans $\mathbb{C}P^{n}$ avec multiplicit\'{e}s
coup\'{e}es et avec "peu de" cibles. On donne aussi un th\'{e}or\`{e}me de 
d\'{e}g\'{e}n\'{e}ration lin\'{e}aire pour des telles applications avec 
multiplicit\'{e}s coup\'{e}es et avec des cibles mobiles. Les preuves utilisent des techniques
de la distribution des valeurs.
\end{abstract}

\section{Introduction}

The uniqueness problem of meromorphic mappings under a condition
on the inverse images of divisors was first studied by R. Nevalinna [8]. He
showed that for two nonconstant meromorphic functions $f$ and $g$ on the
complex plane $\mathbb{C}$, if they have the same inverse images for five
distinct values then $f\equiv g$ . In 1975, H. Fujimoto [3] generalized
Nevanlinna's result to the case of meromorphic mappings of $\mathbb{C}^{m}$
into $\mathbb{C}P^{n}$. He showed  that for two linearly nondegenerate meromorphic
mappings $f$ and $g$ of $\mathbb{C}^{m}$ into $\mathbb{C}P^{n}$, if they
have the same inverse images counted with multiplicities for $(3n+2)$
hyperplanes in general position in $\mathbb{C}P^{n},$ then $f\equiv g$.
Since that time, this problem has been studied intensively by H.Fujimoto
([4], [5] ...), L. Smiley [11], S. Ji [6], \ M. Ru [10], \ Z. Tu [12] and others.

Let $f$ be linearly nondegenerate meromorphic mappings of $\mathbb{C}^{m}$
into $\mathbb{C}P^{n}$. For each hyperplane $H$ we denote by $v_{(f,H)}$ the
map of $\mathbb{C}^{m}$ into $\mathbb{N}_{0}$ such that $v_{(f,H)}(a)$ ($%
a\in \mathbb{C}^{m}$ ) is the intersection multiplicity of the image of $f$
and $H$ at $f(a)$ .

Take $q$ hyperplanes $H_{1},...,H_{q}$ in $\mathbb{C}P^{n}$ in general
position and a positive integer $l_{0}$ .\\

We consider the family $F(\{H_{j}\}_{j=1}^{q},f,l_{0})$ of all linearly
nondegenerate meromorphic mappings $g:$ $\mathbb{C}^{m}$ $\longrightarrow $ $%
\mathbb{C}P^{n}$ satisfying the conditions:\\

(a) $\min \big\{v_{(g,H_{j})},l_{0}\big\}=\min \big\{v_{(f,H_{j})},l_{0}%
\big\}$ for all $j\in \{1,...,q\}$ ,

(b) $\dim \left( f^{-1}(H_{i})\cap f^{-1}(H_{j})\right) \leq m-2$ , for all $%
1\leq i<j\leq q\,,$\ and

(c) $g=f$ on $\bigcup\limits_{j=1}^{q}f^{-1}(H_{j}).$

In 1983, L.Smiley showed that:\\

\noindent \textbf{Theorem A.} ([11])
\textit{If }$q\geq 3n+2$\textit{\ then }$g_{1}$\textit{=}$%
g_{2} $\textit{\ for any }$g_{1}$\textit{,}$g_{2}\in
F(\{H_{j}\}_{j=1}^{q},f,1)$.\\

For the case $q=3n+1$ in [4],[5],[6] the authors gave the following results:\\

\noindent \textbf{Theorem B.} ([6])
\textit{Assume that }$q=3n+1$\textit{.Then for three mappings }$%
g_{1},g_{2},g_{3}\in F(\{H_{j}\}_{j=1}^{q},f,1)$\textit{, the map }$\mathit{g%
}_{1}\mathit{\times g}_{2}\mathit{\times g}_{3}$ $:$ $\mathbb{C}^{m}$ $%
\longrightarrow $ $\mathbb{C}P^{n}\times $ $\mathbb{C}P^{n}\times $ $\mathbb{%
C}P^{n}$ \textit{is algebraically degenerate, namely, }$\{$\textit{\ }$\left(
g_{1}(z),g_{2}(z),g_{3}(z)\right) $\textit{\ }$,z\in C^{m}\}$\textit{\ is
included in a proper algebraic subset of} $\mathbb{C}P^{n}\times $ $\mathbb{C%
}P^{n}\times $ $\mathbb{C}P^{n}.$\\

\noindent \textbf{Theorem C.} ([4])
\textit{Assume that }$q=3n+1$\textit{.Then there are at most two distinct
mappings in }$F(\{H_{j}\}_{j=1}^{q},f,2)$\textit{.\\}

\noindent \textbf{Theorem D.} ([5])
\textit{Assume that }$n=2,q=7.$\textit{\ Then there exist some positive
integer }$l_{0}$\textit{\ and a proper algebraic set }$V$\textit{\ in the
cartesian product of seven copies of the space }$\left( \mathbb{C}%
P^{2}\right) ^{\ast }$\textit{\ of all hyperplanes in }$\mathbb{C}P^{2}$%
\textit{\ such that, for an arbitrary set }$(H_{1},...,H_{7})\notin V$%
\textit{\ and two nondegenerate meromorphic mappings }$f,g$\textit{\ of }$%
\mathbb{C}^{m}$\textit{\ into }$\mathbb{C}P^{2}$\textit{\ with }$\min \big\{%
v_{(g,H_{j})},l_{0}\big\}=\min \big\{v_{(f,H_{j})},l_{0}\big\}$\textit{\ for
all }$j\in \{1,...,7\}$\textit{, we have }$f=g$.\\

In [5], H.Fujimoto also gave some open questions:

+) Does Theorem D remain valid under the assumption that the $H_{j}$ 's are
in general position ?

+) What is a generalization of Theorem D for the case $n \geq 3$ ?\\

In connection with the above results, it is also an interesting problem to ask
whether these results remain valid if the number of hyperplanes is replaced
by a smaller one. In this paper, we will try to get some partial \ answers
to this problem. We give a uniqueness theorem for the case $q\geq n+I\left( 
\sqrt{2n(n+1)}\right) +1$ and a theorem of the linear degeneration for the
case of $(2n+2)$ moving targets (where we denote $I\left( x\right)
:=\min \{k\in $ $\mathbb{N}_{0}:k>x\}$ for a positive constant $x$).

Let $f$, $a$ be two meromorphic mappings of $\mathbb{C}^{m}$ into $\mathbb{C}%
P^{n}$ with reduced representations $f=(f_{0}:\dots :f_{n})$, $%
a=(a_{0}:\dots :a_{n})$.\\
Set $(f,a):=a_{0}f_{0}+\dots +a_{n}f_{n}$. We say that $a$ is ``small'' with
respect to $f$ if $T_{a}(r)= o (T_{f}(r))$ as $r\rightarrow \infty $.
Assuming that $(f,a)%
\not\equiv%
%
0$, we denote by $v_{(f,a)}$ the map of $\mathbb{C}^{m}$ into $\mathbb{N}_{0}$
with $v_{(f,a)}(z)=0$ if $(f,a)(z)\neq 0$ and $v_{(f,a)}(z)=k$ if $z$ is a
zero point of $(f,a)$ with multiplicity $k$.

Let $a_{1},\dots ,a_{q}$ \ $(q\geq n+1)$ be meromorphic mappings of $\mathbb{%
C}^{m}$ into $\mathbb{C}P^{n}$ with reduced representations $%
a_{j}=(a_{j0}:\dots :a_{jn})$, \ $j=1,\dots ,q$. We say that $\big\{a_{j}%
\big\}_{j=1}^{q}$ are in general position if for any $1\leq j_{0}<\dots
<j_{n}\leq q$, \ $\text{det}(a_{j_{k}i},0\leq k,i\leq n)%
\not\equiv%
%
0$. 

For each $j\in\{1,...,q\}$, we put $\widetilde{a_j}=
(\dfrac{a_{j0}}{a_{jt_j}}:...:\dfrac{a_{jn}}{a_{jt_j}})$ and 
$(f,\widetilde{a_j}) =f_0\dfrac{a_{j0}}{a_{jt_j}}+...+f_n\dfrac{a_{jn}}{a_{jt_j}}$
where $a_{jt_j}$ is the first element of $a_{j0}, ..., a_{jn}$ not identically equal to zero.
Let $%
\mathcal{M}%
%
$ be the field (over $\mathbb{C}$ ) of all meromorphic functions on $\mathbb{%
C}^{m}$. Denote by $\mathcal{R}\Big(\big\{a_{j}\big\}_{j=1}^{q}\Big)\subset 
\mathcal{M}%
%
$ the subfield generated by the set $\{\dfrac{a_{ji}}{a_{jt_j}},0\leq i\leq n,1\leq j\leq q\}$
over $\mathbb{C}$ . This subfield is independant of the reduced representations
$a_{j}=(a_{j0}:\dots :a_{jn})$, \ $j=1,\dots ,q$, and it is of course also independant of our
choice of the  ${a_{jt_j}}$, because it contains all quotients of the quotients 
$\dfrac{a_{ji}}{a_{jt_j}}, i=0, \dots ,n$.

We say that $f$ is linearly nondegenerate over $\mathcal{R}\Big(\big\{a_{j}%
\big\}_{j=1}^{q}\Big)$ if $f_{0},\dots ,f_{n}$ are linearly independant over 
$\mathcal{R}\Big(\big\{a_{j}\big\}_{j=1}^{q}\Big)$.

Denote by $\Psi $ the Segre embedding of $\mathbb{C}P^{n}\times $ $\mathbb{C}%
P^{n}$ into $\mathbb{C}P^{n^{2}+2n}$ which is defined by sending the ordered
pair $\left( (w_{0},...,w_{n}),(v_{0},...,v_{n})\right) $ to $%
(...,w_{i}v_{j},...)$ in lexicographic order.

Let $h:\mathbb{C}^{m}\longrightarrow $\textit{\ } $\mathbb{C}P^{n}\times $ $%
\mathbb{C}P^{n}$ be a meromorphic mapping. Let $(h_{0}: \cdots :h_{n^{2}+2n})$ be
a representation of $\Psi  \circ h$ . We say that $h$ is linearly degenerate (with
the algebraic structure in $\mathbb{C}P^{n}\times $ $\mathbb{C}P^{n}$ given
by the Segre embedding) if $h_{0},...,h_{n^{2}+2n}$ are linearly dependant
over $\mathcal{R}\Big(\big\{a_{j}\big\}_{j=1}^{q}\Big).$\\

Our main results are stated as follows:
Let $n,x,y,p$ be nonnegative integers. Assume that: 
\begin{equation*}
2\leq p\leq n\text{ , }1\leq y\leq 2n\text{ , and }
\end{equation*}
\begin{equation*}
0\leq x<\min \{2n-y+1,\frac{(p-1)y}{n+1+y}\}.
\end{equation*}
\\
Let $k$ be an integer or $+ \infty$ with $\frac{2n(n+1+y)(3n+p-x)}{(p-1)y-x(n+1+y)}\leq
k\leq +\infty .$\\

\noindent \textbf{Theorem 1.}
\textit{Let }$f,g$\textit{\ be two linearly nondegenerate meromorphic
mappings of }$\mathbb{C}^{m}$ \textit{into} $\mathbb{C}P^{n}$ \textit{and} $%
\{H_{j}\}_{j=1}^{q}$\textit{\ be }$q:=3n+1-x$\textit{\ hyperplanes in} $%
\mathbb{C}P^{n}$ \textit{in general position}.

\textit{Set }$E_{f}^{j}:=\big\{z\in C^{m}:0\leq v_{(f,H_{j})}(z)\leq k%
\big\}
$\textit{,}$^{\ast }E_{f}^{j}:=\big\{z\in \mathbb{C}^{m}:0<v_{(f,H_{j})}(z)%
\leq k\big\},$\textit{\ and similarily for }$E_{g}^{j}$\textit{, }${}^{\ast
}E_{g}^{j}$\textit{, }$j=1,\dots ,q.$\textit{\\
Assume that }:

(a) $\min \{v_{(f,H_{j})},1\}=\min \{v_{(g,H_{j})},1\}$ \textit{on} $E_{f}^{j}\cap
E_{g}^{j}$ \textit{for all} $j\in \{n+2+y,\dots ,q\}$, \textit{and}

\ \ \ \ $\ \min \{v_{(f,H_{j})},p\}=\min \{v_{(g,H_{j})},p\}$ \textit{on} $%
E_{f}^{j}\cap E_{g}^{j}$ \textit{for all} $j\in \{1,\dots ,n+1+y\}$,

\ \ \ \ 

(b) $\text{dim}\big({}^{\ast }E_{f}^{i}\cap {}^{\ast }E_{f}^{j}\big)\leq m-2$
, dim$\big({}^{\ast }E_{g}^{i}\cap {}^{\ast }E_{g}^{j}\big)\leq m-2$ 
\textit{for all }$1\leq i<j\leq q$,

(c)$f=g$ \textit{on} $\bigcup\limits_{j=1}^{q}\big({}^{\ast }E_{f}^{j}\cap
{}^{\ast }E_{g}^{j}\big)$ .\\
\textit{Then} $f=g.$\\

We state some corollaries of Theorem 1:\\

+) Take $n\geq 2,y=1,p=2,x=0$ and $k\geq n(n+2)(6n+4)$. Then we have:\\

\noindent \textbf{Corollary 1.}
\textit{Let }$f,g$\textit{\ be two linearly nondegenerate meromorphic
mappings of} $\mathbb{C}^{m}$\textit{\ into} $\mathbb{C}P^{n}(n\geq 2)$ 
\textit{and }$\{H_{j}\}_{j=1}^{3n+1}$\textit{\ be hyperplanes in} $\mathbb{C}%
P^{n}$\textit{\ in general position.\\
Assume that:}

(a) $\min \{v_{(f,H_{j})},1\}=\min \{v_{(g,H_{j})},1\}$ \textit{on} $%
E_{f}^{j}\cap E_{g}^{j}$ \textit{for all} $j\in \{n+3,\dots ,3n+1\}$, \textit{and}

\ \ \ \ $\min \{v_{(f,H_{j})},2\}=\min \{v_{(g,H_{j})},2\}$ \textit{on }$%
E_{f}^{j}\cap E_{g}^{j}$ \textit{for all} $j\in \{1,\dots ,n+2\}$,

\ \ \ \ 

(b) $\text{dim}\big({}^{\ast }E_{f}^{i}\cap {}^{\ast }E_{f}^{j}\big)\leq m-2$
, dim$\big({}^{\ast }E_{g}^{i}\cap {}^{\ast }E_{g}^{j}\big)\leq m-2$\textit{%
\ for all }$1\leq i<j\leq 3n+1,$

(c) $f=g$ \textit{on }$\bigcup\limits_{j=1}^{3n+1}\big({}^{\ast
}E_{f}^{j}\cap {}^{\ast }E_{g}^{j}\big)$ .\\
\textit{Then} $f=g.$\\

 Corollary 1 is an improvement of Theorem C. It is also a kind
of generalization of Theorem D to the case where $n\geq 2$ and the
hyperplanes are in general position.\\

+) Take $n\geq 3,y=n+2,p=3,x=1$ and $k=+\infty $. Then we have:\\

\noindent \textbf{Corollary 2.}
\textit{Let }$f,g$\textit{\ be two linearly nondegenerate meromorphic
mappings of} $\mathbb{C}^{m}$\textit{\ into} $\mathbb{C}P^{n}(n\geq 3)$ 
\textit{and }$\{H_{j}\}_{j=1}^{3n}$\textit{\ be hyperplanes in} $\mathbb{C}%
P^{n}$\textit{\ in general position.}\\
\textit{Assume that:}

(a) $\min \{v_{(f,H_{j})},1\}=\min \{v_{(g,H_{j})},1\}$ \textit{for all }$%
j\in \{2n+4,\dots ,3n\}$\textit{,} \textit{and}

\textit{\ \ \ \ \ }$\min \{v_{(f,H_{j})},3\}=\min \{v_{(g,H_{j})},3\}$\textit{%
\ for all }$j\in \{1,\dots ,2n+3\}$\textit{,}

(b) $\text{dim}\big({}f^{-1}(H_{i})\cap f^{-1}(H_{j})\big)\leq m-2$ \textit{%
for all }$1\leq i<j\leq 3n$\textit{, }

(c) $f=g$ on $\bigcup\limits_{j=1}^{3n}f^{-1}(H_{j})$.\\
\textit{Then }$f=g.$\\

+) Take $n\geq 2,y=I\left( \sqrt{2n(n+1)}\right) ,p=n,x=2n-I\left( \sqrt{%
2n(n+1)}\right) ,k=+\infty .$ Then we have:\\

\noindent \textbf{Corollary 3.}
 \textit{\ Let }$f,g$\textit{\ be two linearly nondegenerate
meromorphic mappings of} $\mathbb{C}^{m}$ \textit{into }$\mathbb{C}%
P^{n}\; (n\geq 2)$ \textit{and }$\{H_{j}\}_{j=1}^{n+I\left( \sqrt{2n(n+1)}%
\right) +1}$\textit{\ be hyperplanes in }$\mathbb{C}P^{n}$\textit{\ in
general position.\\
Assume that:}

(a) $\min \{v_{(f,H_{j})},n\}=\min \{v_{(g,H_{j})},n\}$ \textit{for all} $j\in
\{1,\dots ,n+I\left( \sqrt{2n(n+1)}\right)+1\}$, 

(b) $\text{dim}\big({}f^{-1}(H_{i})\cap f^{-1}(H_{j})\big)\leq m-2$ \textit{for all} 
$1\leq i<j\leq n+I\left( \sqrt{2n(n+1)}\right)+1$, 

(c)$f=g$ \textit{on} $\bigcup\limits_{j=1}^{n+I\left( \sqrt{2n(n+1)}\right)
+1}f^{-1}(H_{j})$.\\
\textit{Then} $f=g.$\\

We finally give a result for moving targets:

\noindent \textbf{Theorem 2. }\textit{\ Let }$f,g:\mathbb{C}^{m}\longrightarrow \mathbb{C}P^{n}$ $%
(n\geq 2)$ \textit{be} \textit{two nonconstant meromorphic mappings with
reduced representations }$f=(f_{0}:...:f_{n})$\textit{\ and }$%
g=(g_{0}:...:g_{n}).$\textit{\\
Let }$\{a_{j}\}_{j=1}^{2n+2}$\textit{\ be ''small'' (with respect to }$f$%
\textit{) meromorphic mappings of} $\mathbb{C}^{m}$ \textit{into }$\mathbb{%
C}P^{n}$ \textit{in} \textit{general position with reduced representations }$%
a_{j}=(a_{j0}:...:a_{jn}),$\textit{\ }$j=1,...,2n+2.$\textit{\ Suppose that }%
$(f,a_{j})%
\not\equiv%
%
0$\textit{, }$(g,a_{j})%
\not\equiv%
%
0$\textit{, }$j=1,\dots ,2n+2$.
Take $M$ an integer or $+ \infty$ with $$3n(n+1)\left( 
\begin{array}{c}
2n+2 \\ 
n+1
\end{array}
\right) ^{2}\left[ \left( 
\begin{array}{c}
2n+2 \\ 
n+1
\end{array}
\right) -2 \right] \leq M\leq +\infty .$$
\textit{\\
Assume that:}

(a) $\min \{v_{(f,a_{j})},M\}=\min \{v_{(g,a_{j})},M\}$ \textit{for all} $j\in
\{1,\dots ,2n+2\},$

(b) $\dim \{z\in \mathbb{C}^{m}:{(f,a_{i})}(z)={(f,a_{j})}(z)=0\}\leq
m-2 $ for all $i\neq j,i\in \{1,...,n+4\},j\in \{1,...,2n+2\}$ , 

(c)\textit{There exist }$\gamma _{j}\in \mathcal{R}\Big(\big\{a_{j}\big\}%
_{j=1}^{2n+2}\Big)$ ( $j=1,...,2n+2$ ) \textit{such that} \\
$\gamma _{j}=\frac{a_{j0}f_{0}+...+a_{jn}f_{n}}{a_{j0}g_{0}+...+a_{jn}g_{n}}$
\textit{on }$\left( \bigcup\limits_{i=1}^{n+4}\{z:(f,a_{i})(z)=0\}\right)
\diagdown \{z:(f,a_{j})(z)=0\}.$\\
\textit{Then the mapping} $f\times g:\mathbb{C}^{m}\longrightarrow \mathbb{C}%
P^{n}\times \mathbb{C}P^{n}$ \textit{is linearly degenerate (with the
algebraic} \textit{structure in} $\mathbb{C}P^{n}\times $ $\mathbb{C}P^{n}$ 
\textit{given by the Segre embedding) over }$\mathcal{R}\Big(\big\{a_{j}%
\big\}_{j=1}^{2n+2}\Big).$\\

\noindent \textbf{Remark.} \textit{The condiction (c) is weaker than the following easier one}:

(c')\ \ $f=g$ \textit{on} \ $\bigcup\limits_{i=1}^{n+4}\{z:(f,a_{i})(z)=0\}.$ \\

We finally remark that we also obtained uniqueness theorems with moving targets
(in [1]), and with fixed targets, but  not taking into account, at all, truncations from some
fixed order on (in [2]). But in both cases the number of targets has to be bigger than in
our results above.\\

\noindent \textbf{Acknowledgements.} The second author would like to thank Professor Do Duc
Thai for valuable discussions, the Universit\'{e} de Bretagne Occidentale
(U.B.O.) for its hospitality and for support, and the PICS-CNRS ForMathVietnam
for support.

\section{Preliminaries}

We set $\Vert z\Vert =\big(|z_{1}|^{2}+\dots +|z_{m}|^{2}\big)^{1/2}$ for $%
z=(z_{1},\dots ,z_{m})\in \mathbb{C}^{m}$ and define 
\begin{equation*}
B(r):=\{z\in \mathbb{C}^{m}:|z|<r\},\quad S(r):=\{z\in \mathbb{C}%
^{m}:|z|=r\}\ \text{for all}\ 0<r\leq \infty .
\end{equation*}
Define \ $d^{c}:=\dfrac{\sqrt{-1}}{4\pi }(\overline{\partial }-\partial )$, $%
\ \upsilon :=(dd^{c}\Vert z\Vert ^{2})^{m-1}$ and 
\begin{equation*}
\sigma :=d^{c}\text{log}\Vert z\Vert ^{2}\wedge (dd^{c}\text{log}\Vert
z\Vert ^{2})^{m-1}.
\end{equation*}
Let $F$ be a nonzero holomorphic function on $\mathbb{C}^{m}$. For every $%
a\in \mathbb{C}^{m}$, expanding $F$ as $F=\sum P_{i}(z-a)$ with homogeneous
polynomials $P_{i}$ of degree $i$ around $a$, we define 
\begin{equation*}
v_{F}(a):=\min \{i:P_{i}%
\not\equiv%
%
0\}.
\end{equation*}
Let $\varphi $ be a nonzero meromorphic function on $\mathbb{C}^{m}$. We define the
map $v_{\varphi }$ as follows: for each $z\in \mathbb{C}^{m}$, we choose
nonzero holomorphic functions $F$ and $G$ on a neighborhood $U$ of $z$ such
that $\varphi =\dfrac{F}{G}$ on $U$ and $\text{dim}\big(F^{-1}(0)\cap
G^{-1}(0)\big)\leq m-2$, and then we put $v_{\varphi }(z):=v_{F}(z)$.

Set \ $|v_{\varphi }|:=\overline{\big\{z\in \mathbb{C}^{m}:v_{\varphi
}(z)\neq 0\big\}}$ .

Let $k$, $M$ be positive integers or $+\infty $.

Set 
\begin{equation*}
{}^{\leq M}v_{\varphi }^{[k]}(z)=0\text{ if}\quad v_{\varphi }(z)>M\text{
and }^{\leq M}v_{\varphi }^{[k]}(z)=\min \{v_{\varphi }(z),k\}\text{ }\text{%
if}\quad v_{\varphi }(z)\leq M
\end{equation*}
\begin{equation*}
{}^{>M}v_{\varphi }^{[k]}(z)=0\text{ }\text{if}\quad v_{\varphi }(z)\leq M%
\text{ and }^{>M}v_{\varphi }^{[k]}(z)=\min \{v_{\varphi }(z),k\}\text{ if }%
v_{\varphi }(z)>M.
\end{equation*}
We define 
\begin{equation*}
{}^{\leq M}N_{\varphi }^{[k]}(r):=\int\limits_{1}^{r}\frac{{}^{\leq M}n(t)}{%
t^{2m-1}}dt
\end{equation*}
and 
\begin{equation*}
{}^{>M}N_{\varphi }^{[k]}(r):=\int\limits_{1}^{r}\frac{{}^{>M}n(t)}{t^{2m-1}}%
dt\qquad (1\leq r<+\infty )
\end{equation*}
where, 
\begin{align*}
{}^{\leq M}n(t):={\int\limits_{|v_{\varphi }|\cap B(r)}}{}^{\leq
M}v_{\varphi }^{[k]}.\upsilon \quad & \text{for}\quad m\geq 2{}\text{ , }%
^{\leq M}n(t):={\sum\limits_{|z|\leq t}}{}^{\leq M}v_{\varphi
}^{[k]}(z)\quad \text{for}\quad m=1 \\
{}^{>M}n(t):={\int\limits_{|v_{\varphi }|\cap B(r)}}{}^{>M}v_{\varphi
}^{[k]}.\upsilon \quad & \text{for}\quad m\geq 2\text{ , }{}^{>M}n(t):={%
\sum\limits_{|z|\leq t}}{}^{>M}v_{\varphi }^{[k]}(z)\quad \text{for}\quad
m=1.
\end{align*}

Set \ $N_{\varphi }(r):={}^{\leq \infty }N_{\varphi }^{[\infty ]}(r)$, \ $%
N_{\varphi }^{[k]}(r):={}^{\leq \infty }N_{\varphi }^{[k]}(r)$.

We have the following Jensen's formula (see [5], p.177, observe that his definition of 
$N_{\varphi }(r)$ is a different one than ours): 
\begin{equation*}
N_{\varphi }(r)-N_{\frac{1}{\varphi }}(r)=\int\limits_{S(r)}\text{log}%
|\varphi |\sigma -\int\limits_{S(1)}\text{log}|\varphi |\sigma , \;\;\; 1 \leq r \leq \infty.
\end{equation*}

Let $f : \mathbb{C}^m \longrightarrow \mathbb{C} P^n$ be a meromorphic
mapping. For arbitrary fixed homogeneous coordinates $(w_0 : \dots : w_n)$
of $\mathbb{C} P^n$, we take a reduced representation $f = (f_0 : \dots :
f_n)$ which means that each $f_i$ is a holomorphic function on $\mathbb{C}^m$
and $f(z) = (f_0(z) : \dots : f_n(z))$ outside the analytic set $\{f_0 =
\dots = f_n = 0\}$ of codimension $\geq 2$. Set $\Vert f\Vert = \big(|f_0|^2
+ \dots + |f_n|^2\big)^{1/2}$.

The characteristic function of $f$ is defined by 
\begin{equation*}
T_{f}(r)=\int\limits_{S(r)}\text{log}\Vert f\Vert \sigma -\int\limits_{S(1)}%
\text{log}\Vert f\Vert \sigma ,\quad 1\leq r<+\infty .
\end{equation*}
\\
For a meromorphic function $\varphi $ on $\mathbb{C}^{m}$, the proximity function is defined by
$$m(r,\varphi): = \int\limits_{S(r)}\log^+|\varphi|\sigma$$
 and we have, by the classical First Main Theorem that (see [4], p.135) 
 $$m(r,\varphi) \leq T_{\varphi }(r) +O(1).$$
 Here, the characteristic function $T_{\varphi }(r)$ of $\varphi $ is defined as $\varphi $ can be
considered as a
meroromorphic mapping of $\mathbb{C}^{m}$ into $\mathbb{C}P^{1}.$

We state the First and Second Main Theorem of Value Distribution Theory.
Let $a$ be a meromorphic mapping of $\mathbb{C}^{m}$ into $\mathbb{C}P^{n}$
such that $(f,a)%
\not\equiv%
%
0$,  then for reduced representations   $f=(f_{0}:\dots :f_{n})$ and $%
a=(a_{0}:\dots :a_{n})$, we have:

\vskip0.25cm \noindent \textbf{First Main Theorem.} (Moving target version, see [12], p.569)\textit{
\begin{equation*}
N_{(f,a)}(r)\leq T_{f}(r)+T_{a}(r)+O(1)\quad \text{for}\quad r\geq 1.
\end{equation*}
}For a hyperplane $H:a_{0}w_{0}+\dots +a_{n}w_{n}=0$ in $\mathbb{C}P^{n}$
with $\text{im}\,f\nsubseteq H$, we denote $(f,H)=a_{0}f_{0}+\dots
+a_{n}f_{n}$, where $(f_{0}:\dots :f_{n})$ again is a reduced representation
of $f$.

\vskip0.25cm \noindent \textbf{Second Main Theorem.} (Classical version) \textit{Let $f$ be a
linearly nondegenerate meromorphic mapping of $\mathbb{C}^{m}$ into $\mathbb{%
C}P^{n}$ and $H_{1},\dots ,H_{q}$ $(q\geq n+1)$ hyperplanes of $\mathbb{C}%
P^{n}$ in general position, then 
\begin{equation*}
(q-n-1)T_{f}(r)\leq
\sum\limits_{j=1}^{q}N_{(f,H_{j})}^{[n]}(r)+o(T_{f}(r))\quad
\end{equation*}
}for all$\ r$ \textit{except for a set of finite Lebesgue measure.}

\section{Proof of Theorem 1.}

First of all, we need the following:\\
\noindent%
%

\noindent \textbf{Lemma 3.1.}
\textit{Let }$f,g$\textit{\ be two linearly nondegenerate meromorphic
mappings of} $\mathbb{C}^{m}$ into $\mathbb{C}P^{n}$\textit{\ and }$%
\{H_{j}\}_{j=1}^{q}$\textit{\ be hyperplanes in }$\mathbb{C}P^{n}$\textit{\
in general position. Then there exists a dense subset }$%
\mathcal{C}%
%
$\textit{\ }$\subset $\textit{\ }$\mathbb{C}^{n+1}\diagdown \{0\}$\textit{\
such that for any }$c=(c_{0},...,c_{n})\in 
\mathcal{C}%
%
$\textit{, the hyperplane }$H_{c}$\textit{\ defined by }$c_{0}\omega
_{0}+...+c_{0}\omega _{n}=0$\textit{\ satifies:}
\begin{eqnarray*}
\dim \left( f^{-1}(H_{j})\cap f^{-1}(H_{c})\right) &\leqslant &m-2\mathit{\
\ }\text{\textit{and } }dim\left( g^{-1}(H_{j})\cap g^{-1}(H_{c})\right)
\leqslant m-2\mathit{\ \ } \\
\mathit{\ for\ all\ }j &\in &\{1,...,q\}.
\end{eqnarray*}

\noindent \textbf{Proof.}  We refer to [6], Lemma 5.1. \hfill $\square$\\

We now begin to prove Theorem 1.

Assume that $f%
\not\equiv%
%
g$\textit{\ }.

Let $j_{0}$ be an arbitrarily fixed index, $j_{0}\in \{1,...,n+1+y\}.$ Then
there exists a hyperplane $H$ in $\mathbb{C}P^{n}$ such that: 
\begin{equation*}
\dim \left( f^{-1}(H_{j})\cap f^{-1}(H)\right) \leqslant m-2\mathit{\ \ }%
\text{\textit{, }}\mathit{dim}\left( g^{-1}(H_{j})\cap g^{-1}(H)\right)
\leqslant m-2\mathit{\ \ }
\end{equation*}
\begin{equation}
\mathit{for\ all\ }j\in \{1,...,q\} \mathit{\ and}\text{ }\frac{(f,H_{j_{0}})}{%
(f,H)}%
\not\equiv%
%
\frac{(g,H_{j_{0}})}{(g,H)}:\tag{3.2}
\end{equation}
\\
Indeed, suppose that this assertion does not hold. Then by Lemma 3.1 we have $%
\frac{(f,H_{j_{0}})}{(f,H)}\equiv \frac{(g,H_{j_{0}})}{(g,H)}$ for all
hyperplanes $H$ in $\mathbb{C}P^{n}$. In particular, $\frac{(f,H_{j_{0}})}{%
(f,H_{j_{i}})}\equiv \frac{(g,H_{j_{0}})}{(g,H_{j_{i}})},i\in \{1,...,n\}$
where $\{j_{1},...,j_{n}\}$ is an arbitrary subset of $\{1,...,q\}\diagdown \{j_{0}\}.$
After changing the homogeneous coordinates $(w_{0}:...:w_{n})$ on $\mathbb{C}%
P^{n}$ we may assume that $H_{j_{i}}:w_{i}=0,(i=0,...,n).$ Then $\frac{f_{0}}{%
f_{i}}=\frac{g_{0}}{g_{i}}$ for all $i\in \{1,...,n\}.$ This means that $%
f\equiv g.$ This is a contradiction. Thus we get (3.2).\\
\ \ \ \ \ Since $\min \{v_{(f,H_{j_{0}})},p\}=\min \{v_{(g,H_{j_{0}})},p\}$
on $E_{f}^{j_{0}}\cap E_{g}^{j_{0}}$, $f=g$ on $\bigcup\limits_{j=1}^{q}\big(%
{}^{\ast }E_{f}^{j}\cap {}^{\ast }E_{g}^{j}\big)$ and by (3.2) we have that
a zero point $z_{0}$ of $(f,H_{j_{0}})$ with multiplicity $\leq k$ is either
a zero point of $\frac{(f,H_{j_{0}})}{(f,H)}-\frac{(g,H_{j_{0}})}{(g,H)}$
with multiplicity $\geq \min \{v_{(f,H_{j_{0}})}(z_{0}),p\}$ or a zero point
of $(g,H_{j_{0}})$ with multiplicity $>k$ (outside an analytic set of
codimension $\geq 2)$. \hfill (3.3)

For any $\mathit{\ }j\in \{1,...,q\}\diagdown \{j_{0}\},$ by the asumptions
(a),(c) and by (3.2), we have that a zero point of $(f,H_{j})$ with
multiplicity $\leq k$ is either a zero point of $\frac{(f,H_{j_{0}})}{(f,H)}-%
\frac{(g,H_{j_{0}})}{(g,H)}$ or zero point of $(g,H_{j})$ with multiplicity $%
>k$ (outside an analytic set of codimension $\geq 2).$ \hfill (3.4)

By (3.3) and (3.4), the assumption (b) and by the First Main Theorem we
have
\begin{eqnarray*}
^{\leq k}N_{(f,H_{j_{0}})}^{[p]}+\sum_{j=1,j\neq j_{0}}^{q}{}^{\leq
k}N_{(f,H_{j})}^{[1]}(r) &\leq &N_{\left( \frac{(f,H_{j_{0}})}{(f,H)}-\frac{%
(g,H_{j_{0}})}{(g,H)}\right) }(r)+^{>k}N_{(g,H_{j_{0}})}^{[p]} \\
&&+\sum_{j=1,j\neq j_{0}}^{q}{}^{>k}N_{(g,H_{j})}^{[1]}(r)
\end{eqnarray*}
\begin{equation*}
\leq T_{_{\left( \frac{(f,H_{j_{0}})}{(f,H)}-\frac{(g,H_{j_{0}})}{(g,H)}%
\right) }}(r)+\frac{p}{k+1}N_{(g,H_{j_{0}})}(r)+\frac{1}{k+1}\sum_{j=1,j\neq
j_{0}}^{q}{}N_{(g,H_{j})}(r)+O(1)
\end{equation*}
\begin{equation}
\leq T_{_{\frac{(f,H_{j_{0}})}{(f,H)}}}(r)+T_{\frac{(g,H_{j_{0}})}{(g,H)}%
}(r)+\frac{p+q-1}{k+1}T_{g}(r)+O(1)  \tag{3.5}
\end{equation}
\\

Since dim$(f^{-1}(H_{j_{0}})\cap f^{-1}(H))\leq m-2$ we have:

$
\begin{array}{lll}
\ T_{\frac{(f,H_{j_{0}})}{(f,H)}}(r) & = & \int\limits_{S(r)}\log \,(\left|
(f,H_{j_{0}})\right| ^{2}+\left| (f,H)\right| ^{2})^{\frac{1}{2}}\sigma \ \
+O(1) \\ 
& \leq & \int\limits_{S(r)}\log \,\Vert f\Vert \sigma +O(1)=T_{f}(r)+O(1)\ .
\end{array}
$\\
Similarly,

$
\begin{array}{ccc}
T_{\frac{(g,H_{j_{0}})}{(g,H)}}(r) & \leq & T_{g}(r)+O(1)\ .
\end{array}
$\\
So by (3.5) we have 
\begin{eqnarray*}
^{\leq k}N_{(f,H_{j_{0}})}^{[p]}(r)+\sum_{j=1,j\neq j_{0}}^{q}{}^{\leq
k}N_{(f,H_{j})}^{[1]}(r) &\leq &T_{f}(r)+T_{g}(r) \\
&&+\frac{p+q-1}{k+1}T_{g}(r)+O(1).
\end{eqnarray*}
\\
Similarly, 
\begin{eqnarray*}
^{\leq k}N_{(g,H_{j_{0}})}^{[p]}(r)+\sum_{j=1,j\neq j_{0}}^{q}{}^{\leq
k}N_{(g,H_{j})}^{[1]}(r) &\leq &T_{f}(r)+T_{g}(r) \\
&&+\frac{p+q-1}{k+1}T_{f}(r)+O(1).
\end{eqnarray*}
\\
Thus, 
\begin{eqnarray*}
&&^{\leq k}N_{(f,H_{j_{0}})}^{[p]}(r)+^{\leq
k}N_{(g,H_{j_{0}})}^{[p]}(r)+\sum_{j=1,j\neq j_{0}}^{q}{}\left( ^{\leq
k}N_{(f,H_{j})}^{[1]}(r)+^{\leq k}N_{(g,H_{j})}^{[1]}(r)\right) \\
&\leq &\left( 2+\frac{p+q-1}{k+1}\right) \left( T_{f}(r)+T_{g}(r)\right)
+O(1).
\end{eqnarray*}

\begin{eqnarray*}
&\Rightarrow &\frac{p}{n}\left( ^{\leq k}N_{(f,H_{j_{0}})}^{[n]}(r)+^{\leq
k}N_{(g,H_{j_{0}})}^{[n]}(r)\right) +\frac{1}{n}\sum_{j=1,j\neq
j_{0}}^{q}{}\left( ^{\leq k}N_{(f,H_{j})}^{[n]}(r)+^{\leq
k}N_{(g,H_{j})}^{[n]}(r)\right) \\
&\leq &\frac{2(k+1)+(p+q-1)}{k+1}\left( T_{f}(r)+T_{g}(r)\right) +O(1),
\end{eqnarray*}

(note that $p\leq n).$%
\begin{equation*}
\Rightarrow \frac{p-1}{n}\left( ^{\leq k}N_{(f,H_{j_{0}})}^{[n]}(r)+^{\leq
k}N_{(g,H_{j_{0}})}^{[n]}(r)\right) \leq \frac{2(k+1)+(p+q-1)}{k}\left(
T_{f}(r)+T_{g}(r)\right)
\end{equation*}
\begin{equation}
-\frac{1}{n}\sum_{j=1}^{q}{}\left( ^{\leq k}N_{(f,H_{j})}^{[n]}(r)+^{\leq
k}N_{(g,H_{j})}^{[n]}(r)\right) +O(1)  \tag{3.6}
\end{equation}
\\
By the First and the Second Main Theorem, we have: 
\begin{equation*}
(q-n-1)T_{f}(r)\leq \sum\limits_{j=1}^{q}N_{(f,H_{j})}^{[n]}(r)+o(T_{f}(r))
\end{equation*}
\begin{eqnarray*}
&=&\frac{k}{k+1}\sum_{j=1}^{q}{}^{\leq
k}N_{(f,H_{j})}^{[n]}(r)+\sum_{j=1}^{q}\left( \frac{1}{k+1}\text{ }^{\leq
k}N_{(f,H_{j})}^{[n]}(r)+^{>k}N_{(f,H_{j})}^{[n]}(r)\right) +o(T_{f}(r)) \\
&\leq &\frac{k}{k+1}\sum_{j=1}^{q}{}^{\leq k}N_{(f,H_{j})}^{[n]}(r)+\frac{n}{%
k+1}\sum_{j=1}^{q}{}N_{(f,H_{j})}(r)+o(T_{f}(r)) \\
&\leq &\frac{k}{k+1}\sum_{j=1}^{q}{}^{\leq k}N_{(f,H_{j})}^{[n]}(r)+\frac{nq%
}{k+1}T_{f}(r)+o(T_{f}(r))
\end{eqnarray*}
\begin{equation*}
\Rightarrow \sum_{j=1}^{q}{}^{\leq k}N_{(f,H_{j})}^{[n]}(r)\geq \frac{%
(q-n-1)(k+1)-qn}{k}T_{f}(r)+o(T_{f}(r))
\end{equation*}
\\
Similarly, 
\begin{equation*}
\sum_{j=1}^{q}{}^{\leq k}N_{(g,H_{j})}^{[n]}(r)\geq \frac{(q-n-1)(k+1)-qn}{k}%
T_{g}(r)+o(T_{g}(r))
\end{equation*}
\\
So, 
\begin{equation*}
\sum_{j=1}^{q}{}\left( ^{\leq k}N_{(f,H_{j})}^{[n]}(r)+^{\leq
k}N_{(g,H_{j})}^{[n]}(r)\right) \geq \frac{(q-n-1)(k+1)-qn}{k}\left(
T_{f}(r)+T_{g}(r)\right)
\end{equation*}
\begin{equation}
+o\left( T_{f}(r)+T_{g}(r)\right)  \tag{3.7}
\end{equation}
\\
By (3.6) and (3.7) we have 
\begin{eqnarray*}
&&\frac{p-1}{n}\left( ^{\leq k}N_{(f,H_{j_{0}})}^{[n]}(r)+^{\leq
k}N_{(g,H_{j_{0}})}^{[n]}(r)\right) +o\left( T_{f}(r)+T_{g}(r)\right) \\
&\leq &\left( \frac{2(k+1)+(p+q-1)}{k}-\frac{(q-n-1)(k+1)-qn}{nk}\right)
\left( T_{f}(r)+T_{g}(r)\right)
\end{eqnarray*}
\begin{eqnarray*}
&\Rightarrow &\left( ^{\leq k}N_{(f,H_{j_{0}})}^{[n]}(r)+^{\leq
k}N_{(g,H_{j_{0}})}^{[n]}(r)\right) +o\left( T_{f}(r)+T_{g}(r)\right) \\
&\leq &\frac{(3n+1-q)(k+1)+(2q+p-1)n}{k(p-1)}\left( T_{f}(r)+T_{g}(r)\right) 
\text{ for all }j_{0}\in \{1,...,n+1+y\}
\end{eqnarray*}
\\
So, 
\begin{equation*}
\sum\limits_{j=1}^{n+1+y}\left( ^{\leq k}N_{(f,H_{j})}^{[n]}(r)+^{\leq
k}N_{(g,H_{j})}^{[n]}(r)\right) +o\left( T_{f}(r)+T_{g}(r)\right)
\end{equation*}
\begin{equation}
\leq \frac{(n+1+y)\left[ (3n+1-q)(k+1)+(2q+p-1)n\right] }{k(p-1)}\left(
T_{f}(r)+T_{g}(r)\right)  \tag{3.8}
\end{equation}

By the First and the Second Main Theorem, we have: 
\begin{equation*}
yT_{f}(r)\leq \sum\limits_{j=1}^{n+1+y}N_{(f,H_{j})}^{[n]}(r)+o(T_{f}(r))
\end{equation*}
\begin{eqnarray*}
&=&\frac{k}{k+1}\sum_{j=1}^{n+1+y}{}^{\leq
k}N_{(f,H_{j})}^{[n]}(r)+\sum_{j=1}^{n+1+y}{}\left( \frac{1}{k+1}\text{ }%
^{\leq k}N_{(f,H_{j})}^{[n]}(r)+^{>k}N_{(f,H_{j})}^{[n]}(r)\right)
+o(T_{f}(r)) \\
&\leq &\frac{k}{k+1}\sum_{j=1}^{n+1+y}{}^{\leq k}N_{(f,H_{j})}^{[n]}(r)+%
\frac{n}{k+1}\sum_{j=1}^{n+1+y}{}N_{(f,H_{j})}(r)+o(T_{f}(r)) \\
&\leq &\frac{k}{k+1}\sum_{j=1}^{n+1+y}{}^{\leq k}N_{(f,H_{j})}^{[n]}(r)+%
\frac{n(n+1+y)}{k+1}T_{f}(r)+o(T_{f}(r))
\end{eqnarray*}
\begin{equation*}
\Rightarrow \;\;\;\frac{y(k+1)-n(n+1+y)}{k}T_{f}(r)\leq \sum_{j=1}^{n+1+y}{}^{\leq
k}N_{(f,H_{j})}^{[n]}(r)+o(T_{f}(r)).
\end{equation*}
\\
Similarly,

\begin{equation*}
\frac{y(k+1)-n(n+1+y)}{k}T_{g}(r)\leq \sum_{j=1}^{n+1+y}{}^{\leq
k}N_{(g,H_{j})}^{[n]}(r)+o(T_{g}(r)).
\end{equation*}
\\
So, 
\begin{equation*}
\frac{y(k+1)-n(n+1+y)}{k}\left( T_{f}(r)+T_{g}(r)\right) \leq
\sum_{j=1}^{n+1+y}{}\left( ^{\leq k}N_{(f,H_{j})}^{[n]}(r)+^{\leq
k}N_{(g,H_{j})}^{[n]}(r)\right) 
\end{equation*}
\begin{equation}
+o\left( T_{f}(r)+T_{g}(r)\right)   \tag{3.9}
\end{equation}
\\
By (3.8) and (3.9) we have 
\begin{eqnarray*}
&&\frac{y(k+1)-n(n+1+y)}{k}\left( T_{f}(r)+T_{g}(r)\right) +o\left(
T_{f}(r)+T_{g}(r)\right)  \\
&\leq &\frac{(n+1+y)\left[ (3n+1-q)(k+1)+(2q+p-1)n\right] }{k(p-1)}\left(
T_{f}(r)+T_{g}(r)\right) 
\end{eqnarray*}
So, 
\begin{equation*}
(p-1)\left[ y(k+1)-n(n+1+y)\right] \leq (n+1+y)\left[ x(k+1)+(6n+p+1-2x)n%
\right] 
\end{equation*}
\begin{equation*}
\Rightarrow \;\;\; k+1\leq \frac{2n(n+1+y)(3n+p-x)}{(p-1)y-x(n+1+y)}
\end{equation*}

(note that $(p-1)y-x(n+1+y)>0).$ This is a contradiction. Thus, we have $f\equiv g$. \hfill
$\square $\\

\section{Proof of Theorem 2}

Let $%
\mathcal{G}%
%
$ be a torsion free abelian group and $A=(x_{1},...,x_{q})$ be a $q-$tuple
of elements $x_{i\text{ }}$in $%
\mathcal{G}%
%
$. Let $1<s<r\leq q.$ We say that $A$ has the property $\ P_{r,s}$ if any $r$
elements $x_{p_{1}},...,x_{p_{r}}$ in $A$ satisfy the condition that for any
subset $I\subset \{p_{1},...,p_{r}\}$ with $\#I=s,$ there exists a subset $%
J\subset \{p_{1},...,p_{r}\}$, $J\neq I,\#J=s$ such that $\underset{i\in I}{%
\prod }x_{i}=$ $\underset{j\in J}{\prod }x_{j}.$\\

\noindent \textbf{Lemma 4.1.}
\textit{If} $A$ \textit{has the property} $\ P_{r,s}$, \textit{then there exists a subset} $%
\{i_{1},...,i_{q-r+2}\}\subset \{1,...,q\}$ \textit{such that} $%
x_{i_{1}}=...=x_{i_{q-r+2}}.$\\

\noindent \textbf{Proof.} We refer to [3], Lemma 2.6.\hfill $\square$\\

\noindent \textbf{Lemma 4.2.}
{\it Let $f : \C^m \longrightarrow \C P^n$ be a nonconstant
meromorphic mapping and $\big\{a_i\big\}_{i=0}^n$ be
``small" (with respect to $f$) meromorphic mappings of
$\C^m$ into $\C P^n$ in general position.

Denote the meromorphic mapping,
\begin{align*}
F = \big( c_0 \cdot (f,\tilde a_0) : \cdots :
c_n \cdot (f,\tilde a_n)\big) : \C^m \longrightarrow \C P^n
\end{align*}
where $\big\{c_i\big\}_{i=0}^n$ are ``small" (with respect to $f$)
nonzero meromorphic functions on $\C^m$.

Then,
\begin{align*}
T_F(r) = T_f(r) + o(T_f(r)).
\end{align*}
Moreover, if 
\begin{align*}
f &= (f_0: \cdots : f_n),\\
a_i &= (a_{i0} : \cdots : a_{in}),\\
F &= \Big( \frac{c_0 \cdot (f,\tilde a_0)}{h} : \cdots :
\frac{c_n \cdot (f,\tilde a_n)}{h}\Big)
\end{align*}
are reduced representations, where $h$ is a meromorphic function
on $\C^m$}, then 
\begin{align*}
N_{h}(r) \leq o(T_f(r))
\end{align*}
and
\begin{align*}
N_{\frac{1}{h}}(r) \leq o(T_f(r)).
\end{align*}

\begin{proof}[{\bf Proof}]

Set
\begin{align*}
F_i = \frac{c_i \cdot (f,\tilde a_i)}{h}\,,\quad (i = 0,\dots,n).
\end{align*}
We have
\begin{align}
\begin{cases}
a_{00}f_0 + \dots + a_{0n}f_n = \frac{h}{c_0}F_0 a_{0t_0}\\
\dots \quad \dots \quad \dots \quad \dots \quad \dots \\
a_{n0}f_0 + \dots + a_{nn}f_n =
\frac{h}{c_n} F_n a_{nt_n}\end{cases} \tag{4.3}
\end{align}
Since $(F_0 : \cdots : F_n)$ is a reduced representation of $F$, we have
\begin{align*}
N_{\frac{1}{h}}(r) \leq \sum_{i=0}^n N_{a_{it_i}}(r) +
\sum_{i=0}^n N_{\frac{1}{c_i}}(r) = o(T_f(r)).
\end{align*}
Set
\begin{align*}
P = \begin{pmatrix}
a_{00} & \dots & a_{0n}\cr
\vdots & \ddots & \vdots \cr
a_{n0} & \dots & a_{nn} \end{pmatrix}
\end{align*}
and matrices $P_i$ $(i \in \{0, \dots, n\})$ which are defined from
$P$ after changing the $(i+1)^{th}$ column by
$\begin{pmatrix} F_0 \dfrac{a_{0t_0}}{c_0} \cr
\vdots \cr
F_n \dfrac{a_{nt_n}}{c_n} \end{pmatrix}$.

Put $u_i = \text{det}(P_i)$, $u = \text{det}(P)$, then $u$ is a 
nonzero holomorphic function on $\C^n$ and
\begin{align*}
N_u(r) &= o(T_f(r)),\\
N_{\frac{1}{u_i}}(r) &\leq \sum_{j=0}^n  N_{c_j}(r) = o(T_f(r)),\quad
i = 1,\dots , n.
\end{align*}

By (4.3) we have,
\begin{align}
\begin{cases}
f_0 &= \dfrac{h \cdot u_0}{u} \\
\quad &\vdots \quad \\
f_n &= \dfrac{h \cdot u_n}{u} \end{cases} \tag{4.4}
\end{align}
On the other hand $(f_0 : \cdots : f_n)$ is a reduced 
representation of $f$.

Hence,
\begin{align*}
N_h(r) \leq N_u(r) + \sum_{i=0}^n 
N_{\frac{1}{u_i}}(r) = o(T_f(r)).
\end{align*}
We have

\begin{align}
T_F(r) &= \int\limits_{S(r)} \text{log} 
\Big( \sum_{i=0}^n |F_i|^2\Big)^{1/2}\sigma + 0(1)\notag\\
&= \int\limits_{S(r)} \text{log} 
\Big(\sum_{i=0}^n \Big|\frac{c_i(f,\tilde a_i)}{h}\Big|^2
\Big)^{1/2} \sigma + 0(1)\notag\\
&= \int\limits_{S(r)} \text{log} 
\Big(\sum_{i=0}^n |c_i(f,\tilde a_i)|^2\Big)^{1/2} \sigma -
\int\limits_{S(r)} \text{log} |h| \sigma + 0(1)\notag\\
&\leq \int\limits_{S(r)} \text{log} \Vert f \Vert \sigma +
\int\limits_{S(r)} \text{log} 
\Big(\sum_{i=0}^n |c_i|^2 \cdot \Vert \tilde a_i \Vert^2
\Big)^{1/2} \sigma \notag\\
&\quad  - N_h(r) + N_{\frac{1}{h}}(r) + 0(1) \notag \\
&\leq T_f(r) + \frac{1}{2}\int\limits_{S(r)} \text{log}^+
\Big(\sum_{i=0}^n \Big(
\Big|c_i \frac{a_{i0}}{a_{it_i}}\Big|^2 + \cdots +
\Big| c_i \frac{a_{in}}{a_{it_i}}\Big|^2 \Big)\Big)\sigma 
+o(T_f(r))\notag\\
&\leq T_f(r) + \sum_{i,j=0}^n m \Big(r, c_i \frac{a_{ij}}{a_{it_i}}\Big)
+ o(T_f(r)) \notag\\
&= T_f(r) + o(T_f(r)). \tag{4.5}
\end{align}
(4.4) can be written as
\begin{align*}
\begin{cases}
f_0 &= h \cdot \sum\limits_{i=0}^n b_{i0} F_i\\
\dots &\quad \dots \quad \dots \\
f_n &= h \cdot \sum\limits_{i=0}^n b_{in} F_i \end{cases}
\end{align*}
where $\big\{ b_{ij}\big\}_{i,j=0}^n$ are ``small" (with respect to $f$)
meromorphic functions on $\C^m$.

So,
\begin{align}
T_f(r) &= \int\limits_{S(r)} \text{log} \Vert f \Vert \sigma
+ 0(1) \notag\\
&= \int\limits_{S(r)} \text{log} 
\Big(\sum_{j=0}^n \Big|\sum_{i=0}^n b_{ij}F_i\Big|^2\Big)^{1/2}\sigma +
\int\limits_{S(r)} \text{log} |h| \sigma + 0(1) \notag\\
&\leq \int\limits_{S(r)} \text{log} \Vert F\Vert \sigma +
\int\limits_{S(r)} \text{log}  
\Big( \sum_{i,j} |b_{ij}|^2\Big)^{1/2} \sigma
+ N_h(r) - N_{\frac{1}{h}}(r) + 0(1) \notag \\
&\leq T_F(r) + \int\limits_{S(r)} \text{log}^+
\Big(\sum_{i,j}|b_{ij}|^2\Big)^{1/2} \sigma + o(T_f(r)) \notag\\
&\leq T_F(r) + \sum_{i,j}m(r,b_{ij}) + o(T_f(r)) \notag\\
&= T_F(r) + o(T_f(r)) \tag{4.6}
\end{align}
By (4.5) and (4.6), we get Lemma 4.2.
\end{proof}

We now begin to prove Theorem 2.

The assertion of Theorem 2 is trivial if $f$ or $g$ is linearly degenerate
over $\mathcal{R}\Big(\big\{a_{j}\big\}_{j=1}^{2n+2}\Big).$ So from now we
assume that $f$ and $g$ are linearly nondegenerate over $\mathcal{R}\Big(%
\big\{a_{j}\big\}_{j=1}^{2n+2}\Big).$

Define functions 
\begin{equation*}
h_{j}:=\frac{(a_{j0}f_{0}+...+a_{jn}f_{n})}{(a_{j0}g_{0}+...+a_{jn}g_{n})}%
\,,\quad j\in \{1,\dots ,2n+2\}.
\end{equation*}
For each subset $I\subset \{1,\dots ,2n+2\}$, $\#I=n+1$, set $h_{I}=\prod\limits_{i%
\in I}h_{i},$ $\gamma _{I}=\prod\limits_{i\in I}\gamma _{i}.$\\
Let $%
\mathcal{M}%
%
^{\ast }$ be the abelian multiplication group of all nonzero meromorphic
functions on $\mathbb{C}^{m}$. Denote by $%
\mathcal{H}%
%
$ $\subset 
\mathcal{M}%
%
^{\ast }$ the set of all $h\in 
\mathcal{M}%
%
^{\ast }$ with $h^{k}\in \mathcal{R}\Big(\big\{a_{j}\big\}_{j=1}^{2n+2}\Big)$
for some positive integer $k$. It is easy to see that $%
\mathcal{H}%
%
$ is a subgroup of $%
\mathcal{M}%
%
^{\ast }$ and the multiplication group $%
\mathcal{G}%
%
:=%
\mathcal{M}%
%
^{\ast }\diagup 
\mathcal{H}%
%
$ is a torsion free abelian group. We denote by $[h]$ the class in $%
\mathcal{G}%
%
$ containing $h\in 
\mathcal{M}%
%
^{\ast }$.

We now prove that:\\
 $A:=\left( [h_{1}],...,[h_{2n+2}]\right) $ has the property $ \ \ \ \ \ 
P_{2n+2,n+1}\ .$ \hfill (4.7)

We have
\begin{equation*}
\left\{ 
\begin{array}{c}
a_{j0}f_{0}+\dots +a_{jn}f_{n}=h_{j}(a_{j0}g_{0}+\dots +a_{jn}g_{n}) \\ 
j\in \{1,...,2n+2\}
\end{array}
\right.
\end{equation*}
\begin{equation*}
\Rightarrow \left\{ 
\begin{array}{c}
a_{j0}f_{0}+\dots +a_{jn}f_{n}-h_{j}a_{j0}g_{0}-\dots -h_{j}a_{jn}g_{n}=0 \\ 
1\leq j\leq 2n+2
\end{array}
\right.
\end{equation*}
Therefore, 
\begin{equation*}
\text{det}(a_{j0},\dots ,a_{jn},h_{j}a_{j0},\dots ,h_{j}a_{jn},\ 1\leq j\leq
2n+2)\equiv 0\text{ .}
\end{equation*}
For each $I=\{i_{0},\dots ,i_{n}\}\subset \{1,\dots ,2n+2\}$, $1\leq
i_{0}<\dots <i_{n}\leq 2n+2$, we define 
\begin{equation*}
A_{I}=\frac{(-1)^{\frac{n(n+1)}{2}+i_{0}+\dots +i_{n}}\cdot \text{det}%
(a_{i_{r}j},0\leq r,j\leq n)\cdot \text{det}(a_{i_{s}^{\prime }j},0\leq
s,j\leq n)}{a_ {j_1t_{j_1}} \dots a_ {j_{2n+2}t_{j_{2n+2}}}}
\end{equation*}
where $\{i_{0}^{\prime },\dots ,i_{n}^{\prime }\}=\{1,\dots ,2n+2\}\setminus
\{i_{0},\dots ,i_{n}\}$, $i_{0}^{\prime }<\dots <i_{n}^{\prime }$. We have $%
A_{I}\in \mathcal{R}\Big(\big\{a_{j}\big\}_{j=1}^{2n+2}\Big)$ and $A_{I}%
\not\equiv%
%
0$ , since $\{a_{j}\}_{j=1}^{2n+2}$ are in general position.

Set $L=\{I\subset \{1,\dots ,2n+2\},\#I=n+1\}$, then $\#L=\left( 
\begin{array}{c}
2n+2 \\ 
n+1
\end{array}
\right) $.\\
By the Laplace expansion Theorem, we have 
\begin{equation}
\sum_{I\in L}A_{I}h_{I}\equiv 0.  \tag{4.8}
\end{equation}
\\
\ \ We introduce an equivalence relation on $L$ as follows: $I\backsim J$ if
and only if $\frac{h_{I}}{h_{J}}\in \mathcal{R}\Big(\big\{a_{j}\big\}%
_{j=1}^{2n+2}\Big).$\\
Set $\{L_{1},...,L_{s}\}=L\diagup _{\backsim }$ , ( $s\leq \left( 
\begin{array}{c}
2n+2 \\ 
n+1
\end{array}
\right) $).\\
For each $v\in \{1,\dots ,s\}$, choose $I_{v}\in L_{v}$ and set 
\begin{equation*}
\sum_{I\in L_{v}}A_{I}h_{I}=B_{v}h_{I_{v}},\quad B_{v}\in \mathcal{R}\Big(%
\big\{a_{j}\big\}_{j=1}^{2n+2}\Big).
\end{equation*}
Then (4.8) can be written as 
\begin{equation}
\sum_{v=1}^{s}B_{v}h_{I_{v}}\equiv 0.  \tag{4.9}
\end{equation}
\\
\textbf{Case 1.} There exists some $B_{v}%
\not\equiv%
%
0.$ Without loss of generality we may assume that $B_{v}%
\not\equiv%
%
0$ , for all $v\in \{1,\dots ,l\},B_{v}\equiv 0$ for all $v\in \{l+1,\dots
,s\},(1\leq l\leq s).$\\
By (4.9) we have
\begin{equation}
\sum_{v=1}^{l}B_{v}h_{I_{v}}\equiv 0.  \tag{4.10}
\end{equation}
\\
Denote by $P$ the set of all positive integers $p\leq l$ such that there
exist a subset $K_{p}\subseteq \{1,...,l\},\#K_{p}=p$ and nonzero constants 
$\{c_{i}\}_{i\in K_{p}}$ with $\underset{i\in K_{p}}{\sum }%
c_{i}B_{i}h_{I_{i}}\equiv 0.$ It is clear that $l\in P$ by (4.10). Let $t$ be
the smallest integer in $P,(t\leq l\leq \left( 
\begin{array}{c}
2n+2 \\ 
n+1
\end{array}
\right) ).$\\
We may assume that $K_{t}=\{1,...,t\}.$ Then there exist nonzero constants $%
c_{v},(v=1,...,t)$ such that 
\begin{equation}
\sum_{v=1}^{t}c_{v}B_{v}h_{I_{v}}\equiv 0.  \tag{4.11}
\end{equation}
\\
Since $\frac{h_{I_{i}}}{h_{I_{j}}}\notin \mathcal{R}\Big(\big\{a_{j}\big\}%
_{j=1}^{2n+2}\Big)$ and $h_{I_{i}}%
\not\equiv%
%
0$ for all $1\leq i\neq j\leq t,$ we have $t\geq 3.$\\

Set $\varphi _{1}:=(B_{1}h_{I_{1}}:...:B_{t-1}h_{I_{t-1}})$ , $\varphi
_{2}:=(B_{2}h_{I_{2}}:...:B_{t}h_{I_{t}})$, $\varphi
_{3}:=(B_{1}h_{I_{1}}:B_{3}h_{I_{3}}:...:B_{t}h_{I_{t}})$. They are
meromorphic mappings of $\mathbb{C}^{m}$ into $\mathbb{C}P^{t-2}.$\\
Since $t=\min P$ , we have that $\varphi _{1},\varphi _{2},\varphi _{3}$ are
linearly nondegenerate (over $\mathbb{C)}$ .\\
Without loss of generality, we may assume that 
\begin{equation*}
T_{\varphi _{1}}(r)=\max \{T_{\varphi _{1}}(r),T_{\varphi
_{2}}(r),T_{\varphi _{3}}(r)\}\text{ for all }r\in E,
\end{equation*}
where $E$ is a subset of $[1,+\infty )$ with infinite Lesbesgue measure.%
\\
Since $t\geq 3$ and by the First Main Theorem, we have 
\begin{eqnarray*}
T_{\varphi _{1}}(r) &\geq &\frac{1}{3}\left( T_{\varphi _{1}}(r)+T_{\varphi
_{2}}(r)+T_{\varphi _{3}}(r)\right) \\
&\geq &\frac{1}{3}\left( T_{\frac{B_{1}h_{I_{1}}}{B_{2}h_{I_{2}}}}(r)+T_{%
\frac{B_{2}h_{I_{2}}}{B_{3}h_{I_{3}}}}(r)+T_{\frac{B_{3}h_{I_{3}}}{%
B_{1}h_{I_{1}}}}(r)\right) \\
&\geq &\frac{1}{3}\left( T_{\frac{h_{I_{1}}}{h_{I_{2}}}}(r)+T_{\frac{%
h_{I_{2}}}{h_{I_{3}}}}(r)+T_{\frac{h_{I_{3}}}{h_{I_{1}}}}(r)\right) -\frac{1%
}{3}\left( T_{\frac{B_{1}}{B_{2}}}(r)+T_{\frac{B_{2}}{B_{3}}}(r)+T_{\frac{%
B_{3}}{B_{1}}}(r)\right)
\end{eqnarray*}
\begin{equation}
\geq \frac{1}{3}\left( N_{\frac{h_{I_{1}}}{h_{I_{2}}}-\frac{\gamma _{I_{1}}}{%
\gamma _{I_{2}}}}(r)+N_{\frac{h_{I_{2}}}{h_{I_{3}}}-\frac{\gamma _{I_{2}}}{%
\gamma _{I_{3}}}}(r)+N_{\frac{h_{I_{3}}}{h_{I_{1}}}-\frac{\gamma _{I_{3}}}{%
\gamma _{I_{1}}}}(r)\right) +o(T_{f}(r)),r\in E  \tag{4.12}
\end{equation}

(note that $\frac{h_{I_{i}}}{h_{I_{j}}}%
\not\equiv%
%
\frac{\gamma _{I_{i}}}{\gamma _{I_{j}}}$ since $\frac{h_{I_{i}}}{h_{I_{j}}}%
\notin \mathcal{R}\Big(\big\{a_{j}\big\}_{j=1}^{2n+2}\Big),1\leq i\neq j\leq
3).$\\
Let $(h_{1}^{\prime }:...:h_{t-1}^{\prime })$ be a reduced representation of 
$\varphi _{1}$. Set $h_{t}^{\prime }=\frac{B_{t}h_{I_{t}}h_{1}^{\prime }}{%
B_{1}h_{I_{1}}}.$\\
By (4.11) we have 
\begin{equation}
\sum_{i=1}^{t}c_{i}h_{i}^{\prime }\equiv 0.  \tag{4.13}
\end{equation}
\\
It is easy to see that a zero of $h_{i}^{\prime }\; (i=1,...,t)$ is a zero or a
pole of some $B_{j}h_{I_{j}}$,\thinspace $j\in \{1,\dots ,t\}.$\\
Thus, 
\begin{eqnarray*}
N_{h_{i}^{\prime }}^{[1]}(r) &\leq &\sum_{j=1}^{t}\left(
N_{B_{j}h_{I_{j}}}^{[1]}(r)+N_{\frac{1}{B_{j}h_{I_{j}}}}^{[1]}(r)\right) \\
&\leq &\sum_{j=1}^{t}\left( N_{h_{I_{j}}}^{[1]}(r)+N_{\frac{1}{h_{I_{j}}}%
}^{[1]}(r)\right) +o\left( T_{f}(r)\right) ,i\in \{1,\dots ,t\}.
\end{eqnarray*}

\begin{eqnarray*}
&\Rightarrow &\sum_{i=1}^{t}N_{h_{i}^{\prime }}^{[t-2]}(r)\leq
(t-2)\sum_{i=1}^{t}N_{h_{i}^{\prime }}^{[1]}(r) \\
&\leq &t(t-2)\sum_{j=1}^{t}\left( N_{h_{I_{j}}}^{[1]}(r)+N_{\frac{1}{%
h_{I_{j}}}}^{[1]}(r)\right) +o\left( T_{f}(r)\right)
\end{eqnarray*}
\\
So, by the Second Main Theorem we have 
\begin{eqnarray*}
T_{\varphi _{1}}(r) &\leq &\sum_{i=1}^{t-1}N_{h_{i}^{\prime
}}^{[t-2]}(r)+N_{\left( c_{1}h_{1}^{\prime }+...+c_{t-1}h_{t-1}^{\prime
}\right) }^{[t-2]}(r)+o\left( T_{f}(r)\right) \\
&&\overset{(4.13)}{=}\sum_{i=1}^{t}N_{h_{i}^{\prime }}^{[t-2]}(r)+o\left(
T_{f}(r)\right)
\end{eqnarray*}
\begin{equation}
\leq t(t-2)\sum_{j=1}^{t}\left( N_{h_{I_{j}}}^{[1]}(r)+N_{\frac{1}{h_{I_{j}}}%
}^{[1]}(r)\right) +o\left( T_{f}(r)\right)  \tag{4.14}
\end{equation}
\\

By (4.12) and (4.14) we have 

\begin{equation*}
N_{\frac{h_{I_{1}}}{h_{I_{2}}}-\frac{\gamma _{I_{1}}}{\gamma _{I_{2}}}%
}(r)+N_{\frac{h_{I_{2}}}{h_{I_{3}}}-\frac{\gamma _{I_{2}}}{\gamma _{I_{3}}}%
}(r)+N_{\frac{h_{I_{3}}}{h_{I_{1}}}-\frac{\gamma _{I_{3}}}{\gamma _{I_{1}}}%
}(r)
\end{equation*}
\begin{equation}
\leq 3t(t-2)\sum_{j=1}^{t}\left( N_{h_{I_{j}}}^{[1]}(r)+N_{\frac{1}{h_{I_{j}}%
}}^{[1]}(r)\right) +o\left( T_{f}(r)\right) ,r\in E  \tag{4.15}
\end{equation}
\\
Since $\min \{v_{(f,a_{i})},M\}=\min \{v_{(g,a_{i})},M\}$ for $i\in \{1,\dots
,2n+2\}$, we have 
\begin{equation*}
\{z\in \mathbb{C}^{m}:h_{I_{j}}(z)=0\ \text{or}\ h_{I_{j}}(z)=\infty
\}\subset \bigcup_{i\in I_{j}}\{z\in \mathbb{C}^{m}:v_{(f,a_{i})}(z)>M\
\},j=1,...,t.
\end{equation*}
\\
Thus, 
\begin{equation*}
N_{h_{I_{j}}}^{[1]}(r)+N_{\frac{1}{h_{I_{j}}}}^{[1]}(r)\leq \sum_{i\in
I_{j}}{}^{>M}N_{(f,a_{i})}^{[1]}(r)\leq \frac{1}{M+1}\sum_{i\in
I_{j}}N_{(f,a_{i})}(r)
\end{equation*}
\begin{equation*}
\leq \frac{n+1}{M+1}T_{f}(r)+O(1),j\in \{1,...,t\}
\end{equation*}

( note that $\#I_{j}=n+1).$%
\begin{equation}
\Longrightarrow \sum_{j=1}^{t}\left( N_{h_{I_{j}}}^{[1]}(r)+N_{\frac{1}{%
h_{I_{j}}}}^{[1]}(r)\right) \leq \frac{(n+1)t}{M+1}T_{f}(r)+O(1)  \tag{4.16}
\end{equation}
\\
By (4.15) and (4.16) we have 
\begin{equation*}
N_{\frac{h_{I_{1}}}{h_{I_{2}}}-\frac{\gamma _{I_{1}}}{\gamma _{I_{2}}}%
}(r)+N_{\frac{h_{I_{2}}}{h_{I_{3}}}-\frac{\gamma _{I_{2}}}{\gamma _{I_{3}}}%
}(r)+N_{\frac{h_{I_{3}}}{h_{I_{1}}}-\frac{\gamma _{I_{3}}}{\gamma _{I_{1}}}%
}(r)
\end{equation*}
\begin{equation}
\leq \frac{3(n+1)t^{2}(t-2)}{M+1}T_{f}(r)+o(T_{f}(r)),r\in E  \tag{4.17}
\end{equation}
\\
For each $1\leq s<v\leq 3$ , set $V_{sv}=\{1,...,n+4\}\diagdown \left(
(I_{s}\cup I_{v})\diagdown (I_{s}\cap I_{v})\right) .$\\
Since $\dim \{z\in \mathbb{C}^{m}:(f,a_{i})(z)=(f,a_{j})(z)=0\}\leq m-2$ for
all $i\neq j$, $i\in \{1,\dots ,n+4\},$ $j\in \{1,\dots ,2n+2\},$ and $%
\gamma _{j}=h_{j}$ on $\left(
\bigcup\limits_{i=1}^{n+4}\{z:(f,a_{i})(z)=0\}\right) \diagdown
\{z:(f,a_{j})(z)=0\},$ we have: 
\begin{equation*}
N_{\frac{h_{I_{1}}}{h_{I_{2}}}-\frac{\gamma _{I_{1}}}{\gamma _{I_{2}}}%
}(r)\geq \sum_{i\in V_{12}}{}N_{(f,a_{i})}^{[1]}(r).
\end{equation*}
\\
Indeed, let $z_{0}$ be an arbitrary zero point of some $(f,a_{i}),i\in
V_{12}.$ By omitting an analytic set of codimension $\geq 2$, we may assume
that $(f,a_{j})(z_{0})\neq 0$ for all $j\in \{1,\dots ,2n+2\}\diagdown \{i\}.
$ In particular, $(f,a_{j})(z_{0})\neq 0$ for all $j\in (I_{1}\cup
I_{2})\diagdown (I_{1}\cap I_{2}).$ So $\gamma _{j}(z_{0})=h_{j}(z_{0})$ for
all $j\in (I_{1}\cup I_{2})\diagdown (I_{1}\cap I_{2}).$ Consequently, $z_{0}
$ is a zero point of $\frac{h_{I_{1}}}{h_{I_{2}}}-\frac{\gamma _{I_{1}}}{%
\gamma _{I_{2}}}.$ Thus, the above assertion holds.\\
Similarly, 
\begin{equation*}
N_{\frac{h_{I_{2}}}{h_{I_{3}}}-\frac{\gamma _{I_{2}}}{\gamma _{I_{3}}}%
}(r)\geq \sum_{i\in V_{23}}{}N_{(f,a_{i})}^{[1]}(r),\;\;N_{\frac{h_{I_{3}}}{%
h_{I_{1}}}-\frac{\gamma _{I_{3}}}{\gamma _{I_{1}}}}(r)\geq \sum_{i\in
V_{13}}{}N_{(f,a_{i})}^{[1]}(r).
\end{equation*}
\\
It is easy to see that: $V_{12}\cup V_{23}\cup V_{13}=\{1,...,n+4\}.$\\
Thus, 
\begin{equation}
N_{\frac{h_{I_{1}}}{h_{I_{2}}}-\frac{\gamma _{I_{1}}}{\gamma _{I_{2}}}%
}(r)+N_{\frac{h_{I_{2}}}{h_{I_{3}}}-\frac{\gamma _{I_{2}}}{\gamma _{I_{3}}}%
}(r)+N_{\frac{h_{I_{3}}}{h_{I_{1}}}-\frac{\gamma _{I_{3}}}{\gamma _{I_{1}}}%
}(r)\geq \sum_{i=1}^{n+4}{}N_{(f,a_{i})}^{[1]}(r)\geq
\sum_{i=1}^{n+2}{}N_{(f,a_{i})}^{[1]}(r)  \tag{4.18}
\end{equation}
\\
By (4.17) and (4.18) we have 
\begin{equation}
\sum_{i=1}^{n+2}{}N_{(f,a_{i})}^{[1]}(r)\leq \frac{3(n+1)t^{2}(t-2)}{(M+1)n}%
T_{f}(r)+o(T_{f}(r)),r\in E  \tag{4.19}
\end{equation}
\\
We now prove that:\\
\begin{equation}
\frac{1}{n}T_{f}(r)\leq \sum_{i=1}^{n+2}{}N_{(f,a_{i})}^{[1]}(r)+o(T_{f}(r)).
\tag{4.20}
\end{equation}
\\
Set  
\begin{equation*}
N_{n+2}:=
\begin{pmatrix}
\dfrac{a_{10}}{a_{1t_1}} & \dots  & \dfrac{a_{(n+1)0}}{a_{(n+1)t_{n+1}}}\cr\vdots  & \ddots  & \vdots \cr  \dfrac{a_{1n}}{a_{1t_1}} & \dots 
& \dfrac{a_{(n+1)n}}{a_{(n+1)t_{n+1}}}
\end{pmatrix}
,
\end{equation*}
\\
and matrices $N_{i}$ $(i\in \{1,\dots ,n+1\})$ which are defined by $%
N_{n+2}$ after changing the $i^{th}$ column by $
\begin{pmatrix}
\dfrac{a_{(n+2)0}}{a_{(n+2)t_{n+2}}}\cr\vdots \cr \dfrac{a_{(n+2)n}}{a_{(n+2)t_{n+2}}}
\end{pmatrix}
$. Put $c_{i}=\text{det}(N_{i})$, $(i=1,...,n+2),$ then $\{c_{i}%
\}_{i=1}^{n+2}$ are nonzero meromorphic functions on $\mathbb{C}^{m}$ and $%
c_{i}\in \mathcal{R}\Big(\big\{a_{j}\big\}_{j=1}^{2n+2}\Big).$\\
It is easy to see that: 
\begin{equation}
\sum_{i=1}^{n+1}c_{i}(f,\widetilde{a_{i}})=c_{n+2}(f,\widetilde{a_{n+2}}).  \tag{4.21}
\end{equation}
\\
Denote\ by $F$ the meromorphic mapping $%
(c_{1}(f,\widetilde{a_{1}}):...:c_{n+1}(f,\widetilde{a_{n+1}})):$ $\mathbb{C}^{m}\rightarrow \mathbb{%
C}P^{n}.$\\
Since $f$ is linearly nondegenerate over $\mathcal{R}\Big(\big\{a_{j}\big\}%
_{j=1}^{2n+2}\Big)$  and since $%
\{a_{j}\}_{j=1}^{2n+2}$ are in general position, we have that $F$ is
linearly nondegenerate (over $\mathbb{C)}$ .\\
By Lemma 4.2 we have  
\begin{equation*}
T_{F}(r)=T_{f}(r)+o\left( T_{f}(r)\right) .
\end{equation*}
Let $\Big(\dfrac{c_1(f,\widetilde{a_1})}{h}: ...: \dfrac{c_{n+1}(f,\widetilde{a_{n+1}})}{h}\Big)$ 
be a reduced representation of F, where $h$ is a meromorphic function on $\mathbb{C}^m$. 
By  Lemma 4.2 we have
$$N_h(r) = o(T_f(r)), N_{\frac{1}{h}}(r) =o(T_f(r)).$$
By the Second Main Theorem, we have: 
\begin{eqnarray*}
T_{f}(r)+o\left( T_{f}(r)\right)  &=
&T_{F}(r)\leq
\sum_{i=1}^{n+1}{}N_{\frac{c_{i}(f,\widetilde{a_{i}})}{h}}^{[n]}(r)+
N_{\underset{i=1}{\overset{n+1%
}{\sum }}\frac{c_{i}(f,\widetilde{a_{i}})}{h}}^{[n]}(r)+o\left( T_{F}(r)\right)  \\
&&\overset{(4.21)}{=}\sum_{i=1}^{n+2}{}N_{c_{i}(f,\widetilde{a_{i}})}^{[n]}(r)+N_{\frac{1}{h}}(r)+o\left(
T_{F}(r)\right)  \\
&\leq &\sum_{i=1}^{n+2}{}N_{(f,a_{i})}^{[n]}(r)+\sum_{i=1}^{n+2}N_{\frac{1}{a_{it_i}}}(r)+\sum_{i=1}^{n+2}N_{c_i}(r)+o(T_F(r)) \\
&\leq &n\sum_{i=1}^{n+2}{}N_{(f,a_{i})}^{[1]}(r)+o\left( T_{f}(r)\right) .
\end{eqnarray*}
\\
We get (4.20).\\
By (4.19) and (4.20) we have : 
\begin{equation*}
T_{f}(r)\leq \frac{3(n+1)t^{2}(t-2)}{(M+1)}T_{f}(r)+o(T_{f}(r)),r\in E.
\end{equation*}
\\
This contradicts to 
$$M\geq 3n(n+1)\left( 
\begin{array}{c}
2n+2 \\ 
n+1
\end{array}
\right) ^{2}\left[ \left( 
\begin{array}{c}
2n+2 \\ 
n+1
\end{array}
\right) -2\right] \geq 3(n+1)t^{2}(t-2).$$
\textbf{Case 2.} $B_{v}\equiv 0$ for all $v\in \{1,...,s\}.$ Then $\sum_{I\in
L_{v}}A_{I}h_{I}\equiv 0$ for all $v\in \{1,...,s\}.$ On the other hand $%
A_{I}%
\not\equiv%
%
0,h_{I}%
\not\equiv%
%
0.$ Hence, $\#L_{v}\geq 2$ for all $v\in \{1,...,s\}.$\\
So, for each $I\in L,$ there exists $J\in L,J\neq I$ such that $\frac{h_{I}}{%
h_{J}}\in \mathcal{R}\Big(\big\{a_{j}\big\}_{j=1}^{2n+2}\Big)$ . This
implies that $\prod\limits_{i\in I}[h_{i}]=\prod\limits_{i\in J}[h_{j}].$%
\\
We get (4.7). \hfill $\square $\\

By Lemma 4.1 there exist $j_{1},j_{2}\in \{1,...,2n+2\},$ $j_{1}\neq j_{2}$
such that $[h_{j_{1}}]=[h_{j_{2}}].$\\
By the definition, we have $\frac{h_{j_{1}}}{h_{j_{2}}}\in 
\mathcal{H}%
%
.$ This means that $\left( \frac{h_{j_{1}}}{h_{j_{2}}}\right) ^{k}\in 
\mathcal{R}\Big(\big\{a_{j}\big\}_{j=1}^{2n+2}\Big)$ for some positive
integer $k.$\\
So $\left( \frac{(f,a_{j_{1}})(g,a_{j_{2}})}{(g,a_{j_{1}})(f,a_{j_{2}})}%
\right) ^{k}\in \mathcal{R}\Big(\big\{a_{j}\big\}_{j=1}^{2n+2}\Big).$\\
Take $\{i_{1},...,i_{n+2}\}\subseteq \{1,...,n+4\}\diagdown \{j_{1},j_{2}\}.$%
\\
Similarly to (4.20), we have: 
\begin{equation}
\frac{1}{n}T_{f}(r)\leq
\sum_{s=1}^{n+2}{}N_{(f,a_{i_{s}})}^{[1]}(r)+o(T_{f}(r)).  \tag{4.22}
\end{equation}
\\
\ \ +) If $\left( \frac{(f,a_{j_{1}})(g,a_{j_{2}})}{(g,a_{j_{1}})(f,a_{j_{2}})%
}\right) ^{k}-\left( \frac{\gamma _{j_{1}}}{\gamma _{j_{2}}}\right) ^{k}%
\not\equiv%
%
0,$ then by the assumptions (b) and (c) we have
\begin{equation}
N_{\left( \frac{(f,a_{j_{1}})(g,a_{j_{2}})}{(g,a_{j_{1}})(f,a_{j_{2}})}%
\right) ^{k}-\left( \frac{\gamma _{j_{1}}}{\gamma _{j_{2}}}\right)
^{k}}(r)\geq \sum_{s=1}^{n+2}{}N_{(f,a_{i_{s}})}^{[1]}(r)  \tag{4.23}
\end{equation}
\\
Indeed, let $z_{0}\,$be an arbitrary zero point of some $(f,a_{i_{s}})$ ,
( $1\leq s\leq n+2).$ By omitting an analytic set of codimemsion $\geq 2,$
we may assume that $(f,a_{j_{1}})(z_{0})\neq 0$ , $(f,a_{j_{2}})(z_{0})\neq
0 $ (note that $j_{1},j_{2}\neq i_{s}).$ Then $\gamma _{j_{1}}(z_{0}\,)=%
\frac{(f,a_{j_{1}})}{(g,a_{j_{1}})}(z_{0}\,),\; \gamma _{j_{2}}(z_{0}\,)=\frac{%
(f,a_{j_{2}})}{(g,a_{j_{2}})}(z_{0}\,).$ Thus $z_{0}$ is a zero point of $%
\left( \frac{(f,a_{j_{1}})(g,a_{j_{2}})}{(g,a_{j_{1}})(f,a_{j_{2}})}\right)
^{k}-\left( \frac{\gamma _{j_{1}}}{\gamma _{j_{2}}}\right) ^{k}.$ We get
(4.23).\\
By the First Main Theorem and by (4.22),(4.23) we have: 
\begin{equation*}
T_{\left( \frac{(f,a_{j_{1}})(g,a_{j_{2}})}{(g,a_{j_{1}})(f,a_{j_{2}})}%
\right) ^{k}}(r)+T_{\left( \frac{\gamma _{j_{1}}}{\gamma _{j_{2}}}\right)
^{k}}(r)\geq N_{\left( \frac{(f,a_{j_{1}})(g,a_{j_{2}})}{%
(g,a_{j_{1}})(f,a_{j_{2}})}\right) ^{k}-\left( \frac{\gamma _{j_{1}}}{\gamma
_{j_{2}}}\right) ^{k}}(r)
\end{equation*}
\begin{equation*}
\geq \sum_{s=1}^{n+2}{}N_{(f,a_{i_{s}})}^{[1]}(r)\geq \frac{1}{n}%
T_{f}(r)+o(T_{f}(r)).
\end{equation*}
\\
This is a contradiction, since $\gamma _{j_{1}},\gamma _{j_{2}},\left( \frac{%
(f,a_{j_{1}})(g,a_{j_{2}})}{(g,a_{j_{1}})(f,a_{j_{2}})}\right) ^{k}\in 
\mathcal{R}\Big(\big\{a_{j}\big\}_{j=1}^{2n+2}\Big).$\\
Thus, $\left( \frac{(f,a_{j_{1}})(g,a_{j_{2}})}{(g,a_{j_{1}})(f,a_{j_{2}})}%
\right) ^{k}\equiv \left( \frac{\gamma _{j_{1}}}{\gamma _{j_{2}}}\right)
^{k}.$ So, $\frac{(f,a_{j_{1}})(g,a_{j_{2}})}{(g,a_{j_{1}})(f,a_{j_{2}})}%
\equiv \alpha \frac{\gamma _{j_{1}}}{\gamma _{j_{2}}}$ , where $\alpha $ is
a constant. This implies that $f\times g$ is linearly degenerate over $%
\mathcal{R}\Big(\big\{a_{j}\big\}_{j=1}^{2n+2}\Big).$\\
We have completed proof of Theorem 2. \hfill $\square $

 \noindent Gerd Dethloff, Tran Van Tan \\
Universit\'{e} de Bretagne Occidentale \\
  UFR Sciences et Techniques \\
D\'{e}partement de Math\'{e}matiques \\
6, avenue Le Gorgeu, BP 452 \\
  29275 Brest Cedex, France \\
e-mail: gerd.dethloff@univ-brest.fr\\

\end{document}